\documentclass[11pt]{amsart}

\usepackage{amscd}
\usepackage{amsmath, amssymb}
\usepackage{amsfonts}
\usepackage[colorlinks,linkcolor=blue, pdfstartview=FitH]{hyperref}

\numberwithin{equation}{section}

\theoremstyle{plain}
\newtheorem{thm}{Theorem}[section]
\newtheorem{theorem}[thm]{Theorem}
\newtheorem{lemma}[thm]{Lemma}
\newtheorem{corollary}[thm]{Corollary}
\newtheorem{proposition}[thm]{Proposition}

\theoremstyle{definition}

\newtheorem{remark}[thm]{Remark}

\newtheorem{definition}[thm]{Definition}

\newtheorem{example}[thm]{Example}

\newtheorem{defn-thm}[thm]{Definition-Theorem}

\newcommand{\C}{{\mathbb C}}

\renewcommand{\P}{{\mathbb P}}

\newcommand{\R}{{\mathbb R}}

\renewcommand{\S}{{\mathbb S}}

\newcommand{\qtq}[1]{\quad\mbox{#1}\quad}
\newcommand{\bp}{\bar{\partial}}
\newcommand{\Om}{\Omega}

\newcommand{\ts}{\otimes}

\newcommand{\btheorem}{\begin{theorem}}
\newcommand{\etheorem}{\end{theorem}}
\newcommand{\bproposition}{\begin{proposition}}
\newcommand{\eproposition}{\end{proposition}}
\newcommand{\bdefinition}{\begin{definition}}
\newcommand{\edefinition}{\end{definition}}
\newcommand{\bcorollary}{\begin{corollary}}
\newcommand{\ecorollary}{\end{corollary}}
\newcommand{\bproof}{\begin{proof}}
\newcommand{\eproof}{\end{proof}}
\newcommand{\bremark}{\begin{remark}}
\newcommand{\eremark}{\end{remark}}
\newcommand{\eexample}{\end{example}}
\newcommand{\bexample}{\begin{example}}
\newcommand{\la}{\langle}
\newcommand{\elemma}{\end{lemma}}
\newcommand{\blemma}{\begin{lemma}}
\newcommand{\ra}{\rangle}
\newcommand{\sq}{\sqrt{-1}}

\newcommand{\p}{\partial}

\renewcommand{\bar}{\overline}

\renewcommand{\phi}{\varphi}

\newcommand{\ee}{\end{eqnarray*}}
\newcommand{\be}{\begin{eqnarray*}}

\newcommand{\beq}{\begin{equation}}
\newcommand{\eeq}{\end{equation}}

\newcommand{\bd}{\begin{enumerate}}
\newcommand{\ed}{\end{enumerate}}

\renewcommand{\hat}{\widehat}

\renewcommand{\bf}{\textbf}

\renewcommand{\bf}{\textbf}

\newcommand{\ci}{I_i}

\newcommand{\cbj}{I_{ \bar j}}

\renewcommand{\>}{\rightarrow}

\begin{document}
\title{Hermitian Harmonic maps and non-degenerate curvatures}

\makeatletter
\let\uppercasenonmath\@gobble
\let\MakeUppercase\relax
\let\scshape\relax
\makeatother
\author{Kefeng Liu, Xiaokui Yang}
\date{}

\address{Kefeng Liu, Department of Mathematics, UCLA, Los Angeles, CA}
\email{liu@math.ucla.edu}

\address{Xiaokui Yang, Department of Mathematics, Northwestern
University, Evanston, IL} \email{xkyang@math.northwestern.edu}

\maketitle

\begin{abstract}
In this paper, we study the existence of various harmonic maps from
Hermitian manifolds to K\"ahler, Hermitian and Riemannian manifolds
respectively. By using refined Bochner formulas on Hermitian
(possibly non-K\"ahler) manifolds, we derive  new rigidity results
on Hermitian harmonic maps from compact Hermitian manifolds to
Riemannian manifolds, and we also obtain the complex analyticity of
pluri-harmonic maps from compact complex manifolds to compact
K\"ahler manifolds (and Riemannian manifolds) with non-degenerate
curvatures, which are analogous to several fundamental results in
\cite[Siu]{Siu}, \cite[Jost-Yau]{JY} and \cite[Sampson]{Sa2}.


\end{abstract}

\setcounter{tocdepth}{1} \tableofcontents



\section{Introduction}

In the seminal work \cite{Siu} of Siu, he proved that
\btheorem[{\cite[Siu]{Siu}}]\label{siu0}Let $f:(M,h)\>(N,g)$ be a
harmonic map between compact K\"ahler manifolds. If $(N,g)$ has
strongly negative curvature and $rank_{\R} df\geq 4$, then $f$ is
holomorphic or anti-holomorphic. \etheorem

There is a natural question, whether one can obtain similar results
when $(M,h)$ is Hermitian but non-K\"ahler. The main difficulty
arises from the torsion of  non-K\"ahler metrics when applying
Bochner formulas (or Siu's $\p\bp$-trick) on Hermitian manifolds. On
the other hand, it is well-known that if the domain manifold $(M,h)$
is non-K\"ahler, there are various different harmonic maps and they
are mutually different( see Section \ref{hd} for more details). In
particular, holomorphic maps or anti-holomorphic maps are not
necessarily harmonic (with respect to the background Riemannian
metrics). The first result along this line was proved by Jost and
Yau in their fundamental work \cite{JY}, where they used ``Hermitian
harmonic map":

\btheorem[{\cite[Jost-Yau]{JY}}]\label{JY111} Let  $(N,g)$ be a
compact K\"ahler manifold, and $(M,h)$ a compact Hermitian manifold
with $\p\bp\omega^{m-2}_h=0$ where $m=\dim_\C M$. Let
 $f:(M,h)\>(N,g)$ be a Hermitian harmonic map. Then $f$ is holomorphic or
anti-holomorphic if $(N,g)$ has strongly negative curvature (in the
sense of Siu) and $rank_{\R} df\geq 4$. \etheorem

\noindent In the proof of Theorem \ref{JY111}, the condition
$\p\bp\omega^{m-2}_h=0$ plays the key role. Now a Hermitian manifold
$(M,h)$ with  $\p\bp\omega_h^{m-2}=0$ is called
\emph{astheno-K\"ahler}.

In this paper, we study various harmonic maps from general Hermitian
 manifolds and also investigate the complex analyticity of Hermitian
harmonic maps and pluri-harmonic maps.  Let $f:(M,h)\>(N,g)$ be a
smooth map between two compact manifolds. If $(M, h)$ is Hermitian,
we can consider the critical points of the partial energies
$$E''(f)=\int_M|\bp f|^2\frac{\omega_h^m}{m!},\ \ \ \ E'(f)=\int_M|\p f|^2\frac{\omega_h^m}{m!}.$$
If the target manifold $(N,g)$ is  a K\"ahler manifold (resp.
Riemannian manifold), the Euler-Lagrange equations of the partial
energies $E''(f)$ and $E'(f)$ are $\bp_E^*\bp f=0$ and $\p_E^*\p
f=0$ respectively where $E$ is the pullback vector bundle
$f^*(T^{1,0}N)$ (resp. $E=f^*(TN)$). They are called $\bp$-harmonic
and $\p$-harmonic maps respectively.
 In general, $\bp$-harmonic maps are not
necessarily $\p$-harmonic and vice versa. In \cite{JY}, Jost and Yau
considered a reduced harmonic map equation $$
-h^{\alpha\bar\beta}\left(\frac{\p^2 f^i}{\p z^\alpha\p\bar
z^\beta}+\Gamma_{jk}^i\frac{\p f^j}{\p \bar z^\beta}\frac{\p f^k}{\p
z^\alpha}\right)=0. $$ Now it is called \emph{Hermitian harmonic
map} (or \emph{pseudo-harmonic map}). The Hermitian harmonic map has
generalized divergence free structures
$$ \left(\bp_E-2\sq\p^*\omega_h\right)^*\left(\bp f\right)=0
\qtq{or} \left(\p_E+2\sq\bp^*\omega_h\right)^*\left(\p f\right)=0.
$$ The classical harmonic maps, $\bp$-harmonic maps, $\p$-harmonic
maps and Hermitian harmonic maps coincide if the domain manifold
$(M,h)$ is K\"ahler.

 In Section \ref{hd}, we clarify and summarize the definitions of  various harmonic maps
  from Hermitian manifolds to K\"ahler manifolds, to Hermitian manifolds and to
 Riemannain manifolds respectively. Their relations are also discussed.

 By using methods developed in \cite[Eells-Sampson]{ES} and  \cite[Jost-Yau]{JY}, we show in
Section \ref{existence} that $\bp$-harmonic maps and $\p$-harmonic
maps always exist if the target manifold $(N,g)$ has non-positive
Riemannian sectional curvature.

In Section \ref{hdd}, we study Hermitian harmonic maps from
Hermitian manifolds to  Riemannian manifolds. At first, we obtain
the following generalization of a fundamental result of
Sampson(\cite{Sa2}), which is also an analogue to Theorem
\ref{JY111}:

\btheorem Let $(M,h)$ be a compact Hermitian manifold with
$\p\bp\omega^{m-2}_h=0$ and $(N,g)$ a  Riemannian manifold. Let
$f:(M,h)\>(N,g)$ be a Hermitian harmonic map, then $rank_{\R}df\leq
2$ if $(N,g)$ has strongly Hermitian-negative curvature. In
particular, if $\dim_{\C}M>1$, there is no Hermitian harmonic
immersion of $M$ into Riemannian  manifolds of constant negative
curvature. \etheorem

\noindent Here, the strongly Hermitian-negative curvatures on
Riemannian manifolds (see Definition \ref{samdef}) are analogous to
Siu's strongly negative curvatures on K\"ahler manifolds. For
example, Riemannian manifolds with negative constant curvatures have
strongly Hermitian-negative curvatures. On the other hand,
  the condition
$\p\bp\omega_h^{m-2}=0$ can be satisfied on a large class of
Hermitian non-K\"ahler manifolds(\cite{Matsuo}, \cite{FT}), for
example, Calabi-Eckmann manifolds $S^{2p+1}\times \S^{2q+1}$ with
$p+q+1=m$.

In Section \ref{pl},  we  consider the complex analyticity of
pluri-harmonic maps from compact \emph{complex manifolds} to compact
K\"ahler manifolds and  Riemannian manifolds respectively. The
following results are also analogous to Theorem \ref{siu0} and
Theorem \ref{JY111}:

\btheorem\label{i1} Let $f:M\>(N,g)$ be a pluri-harmonic map from a
compact complex manifold $M$ to a compact K\"ahler manifold $(N,g)$.
Then it is holomorphic or anti-holomorphic if $(N,g)$ has
non-degenerate curvature and $rank_{\R} df\geq 4$. \etheorem
\noindent Here, ``non-degenerate curvature" (see Definition
\ref{Siudef}) is a generalization of Siu's ``strongly positive
curvature". For example, both manifolds with strongly positive
curvatures and manifolds with strongly negative curvatures have
non-degenerate curvatures. Hence, in particular,

\bcorollary\label{key222} Let $f:M\>(N,g)$ be a pluri-harmonic map
from a compact complex manifold $M$ to a compact K\"ahler manifold
$(N,g)$. Then it is holomorphic or anti-holomorphic if $(N,g)$ has
strongly negative curvature and $rank_{\R} df\geq 4$.

\ecorollary

\noindent We can see from the proof of Theorem \ref{pluri1} that
Corollary \ref{key222} also holds if the target manifold $N$ is the
compact quotient of a bounded symmetric domain and $f$ is a
submersion. The key ingredients in the proofs are some new
observations on refined Bochner formulas on the vector  bundle
$E=f^*(T^{1,0}N)$ on the Hermitian( possibly non-K\"ahler) manifold
$M$.


As similar as Theorem \ref{i1}, we obtain \btheorem\label{i4} Let
$f:M\>(N,g)$ be a pluri-harmonic map from a compact complex manifold
$M$ to a  Riemannian manifold $(N,g)$. If the Riemannian curvature
$R^g$ of $(N,g)$ is Hermitian non-degenerate at some point $p$, then
$rank_{\R} df(p)\leq 2$. \etheorem

As examples, we show
 \bcorollary \bd\item Any pluri-harmonic map from the Calabi-Eckmann manifold
$\S^{2p+1}\times \S^{2q+1}$ to the real space form $N(c)$ is
constant if $p+q\geq 1$.

\item Any pluri-harmonic map from  $\C\P^n$ to the real space form $N(c)$ is
constant if $n\geq 2$.

 \ed\ecorollary

\bf{Acknowledgement}. The second named author  wishes to thank
Valentino Tosatti for his invaluable suggestions and support.

\section{Connections on vector bundles}
\subsection{Connections on vector bundles}
Let $E$ be a Hermitian {complex} vector bundle or a Riemannian real
vector bundle over a compact Hermitian manifold $(M,h)$ and
$\nabla^E$ be a metric connection on $E$. There is a natural
decomposition $\nabla^E=\nabla^{'E}+\nabla{''^E}$ where \beq
\nabla^{'E}:\Gamma(M,E)\>\Om^{1,0}(M,E) \qtq{and}
\nabla^{''E}:\Gamma(M,E)\>\Om^{0,1}(M,E).
  \eeq
Moreover, $\nabla^{'E}$ and $\nabla^{''E}$  induce two differential
operators. The first one is
 $
\p_E:\Om^{p,q}(M,E)\>\Om^{p+1,q}(M,E) $ defined by \beq \p_E(\phi\ts
s)=\left(\p\phi\right)\ts s+(-1)^{p+q}\phi\wedge \nabla^{'E}s \eeq
for any $\phi\in \Om^{p,q}(M)$ and $s\in \Gamma(M,E)$. The operator
$ \bp_E:\Om^{p,q}(M,E)\>\Om^{p+1,q}(M,E) $  is defined similarly.
For any $\phi\in \Om^{p,q}(M)$ and $s\in \Gamma(M,E)$, \beq
\left(\p_E\bp_E+\bp_E\p_E\right)(\phi\ts s)=\phi\wedge
\left(\p_E\bp_E+\bp_E\p_E\right)s.\eeq  The operator
$\p_E\bp_E+\bp_E\p_E$ is represented by the curvature tensor $R^E
\in \Gamma(M, \Lambda^{1,1}T^*M\ts E^*\ts E)$. For any $\phi,
\psi\in \Om^{\bullet,\bullet}(M,E)$, there is a \emph{sesquilinear
pairing} \beq \left\{ \phi, \psi\right\}=\phi^\alpha\wedge \bar
{\psi^\beta} \langle e_\alpha, e_\beta\rangle \eeq  if
$\phi=\phi^\alpha e_\alpha$ and $\psi=\psi^\beta e_\beta$ in the
local frames $\{e_\alpha\}$ on $E$. By the metric compatible
property of $\nabla^E$, \beq \begin{cases}\p\{\phi,\psi\}=\{\p_E
\phi,\psi\}+(-1)^{p+q}\{\phi, \bp_E\psi\} \\
\bp\{\phi,\psi\}=\{\bp_E \phi,\psi\}+(-1)^{p+q}\{\phi, \p_E\psi\}
\end{cases}\eeq
 if
$\phi\in\Om^{p,q}(M,E)$.
 Let $\omega$ be the fundamental $(1,1)$-form of the Hermitian metric $h$, i.e., \beq
\omega=\frac{\sq}{2}h_{i\bar j}dz^i\wedge d\bar z^j. \eeq
 On the Hermitian manifold $(M,h,\omega)$, the norm on
$\Om^{p,q}(M,E)$ is defined by \beq (\phi,\psi)=\int_M
\{\phi,*\psi\}=\int_M\left(\phi^\alpha\wedge *\bar
{\psi^\beta}\right) \langle e_\alpha, e_\beta\rangle  \eeq for
$\phi,\psi\in \Om^{p,q}(M,E)$. The adjoint operators of
$\p,\bp,\p_E$ and $\bp_E$ are denoted by $\p^*,\bp^*, \p^*_E$ and
$\bp_E^*$ respectively. We shall use the following computational
lemmas frequently in the sequel and  the proofs of them can be found
in \cite{LY12}. \blemma\label{computationlemma}
 We have the following formula:
 \beq\begin{cases} \bp_E^*(\phi\ts s)=(\bp^*\phi)\ts s-h^{i\bar j}\left(\cbj\phi\right) \wedge \nabla_i^Es\\
 \p_E^*(\phi\ts s)=(\p^*\phi)\ts s-h^{i\bar j}\left(\ci\phi\right) \wedge \nabla_{\bar j}^Es \end{cases}\label{bpe*}\eeq
 for any $\phi\in \Om^{p,q}(M)$ and $s\in \Gamma(M,E)$. We use the
 compact notations,
 $$I_{i}=I_{\frac{\p}{\p z^i}},\ \ I_{\bar j}=I_{\frac{\p}{\p\bar z^j}},\ \nabla_i^E=\nabla_{\frac{\p}{\p z^i}}^E,\ \nabla_{\bar j}^E=\nabla^E_{\frac{\p}{\p\bar z^j}}$$
where $I_X$ the contraction operator by the (local) vector field
$X$. \elemma

\blemma\label{realbochner} Let $E$ be a Riemannian real vector
bundle or a Hermitian complex vector bundle over a compact Hermitian
manifold $(M,h,\omega)$. If $\nabla$ is a metric connection on $E$
and $\tau$ is the operator of type $(1,0)$ defined by
$\tau=[\Lambda,2\partial\omega]$ on $\Om^{\bullet,\bullet}(M,E)$,
then we have\bd
\item[(1)] $[\bp_E^*,L]=\sq(\p_E+\tau)$, $[\p^*_E,L]=-\sq(\bp_E+\bar{\tau})$; \item[(2)]
$[\Lambda,\p_E]=\sqrt{-1}(\bp_E^*+\bar{\tau}^{*})$,
 $[\Lambda,\bp_E]=-\sqrt{-1}(\p_E^*+\tau^{*})$. \ed Moreover, \beq
\Delta_{\bp_E}=\Delta_{\p_E}+\sq[R^E,\Lambda]+\left(\tau^*\p_E+\p_E\tau^*\right)-\left(\bar\tau^*\bp_E+\bp_E\bar
\tau^*\right)\label{bochner} \eeq where $\Lambda$ is the contraction
operator by $2\omega$. \elemma

\section{Harmonic map equations }\label{hd}

In this section, we shall  clarify and summarize the definitions of
various harmonic maps between two of the following manifolds:
Riemannian manifolds, Hermitian manifolds and K\"ahler manifolds.
There are many excellent references on this interesting topic, and
we refer the reader to \cite{EL,EL2,EL3} and also the references
therein.

\subsection{Harmonic maps between Riemannian manifolds}

 Let $(M,h)$ and $(N,g)$ be two compact Riemannian
 manifolds. Let $f:(M,h)\>(N,g)$ be a smooth map. If $E=f^*(TN)$,
 then  $df$ can be regarded as  an $E$-valued one form. There
 is an induced connection $\nabla^E$ on $E$ by the Levi-Civita connection on
 $TN$. In the local coordinates $\{x^\alpha\}_{\alpha=1}^m$, $\{y^i\}_{i=1}^n$
 on
 $M$ and $N$ respectively, the local frames of $E$ are denoted by
 $e_i=f^*\left(\frac{\p}{\p y^i}\right)$ and
 $$\nabla^E e_i=f^*\left(\nabla \frac{\p}{\p y^i}\right)=\Gamma_{ij}^k(f)\frac{\p f^j}{\p x^\alpha}dx^\alpha \ts e_k.$$
The connection $\nabla^E$ induces a differential operator $
d_E:\Om^p(M,E)\>\Om^{p+1}(M,E)$ given by $d_E(\phi\ts
s)=\left(d\phi\right)\ts s+(-1)^p\phi\wedge \nabla^E s$ for any
$\phi\in \Om^p(M)$ and $s\in\Gamma(M,E)$.  As a classical result,
the Euler-Lagrange equation of the energy \beq E(f)=\int_M |df|^2
dv_M\eeq is $d_E^* df=0$,  i.e. \beq h^{\alpha\beta}\left(\frac{\p^2
f^i}{\p x^\alpha\p x^\beta}-\frac{\p f^i}{\p
x^\gamma}\Gamma_{\alpha\beta}^\gamma
 +\Gamma_{jk}^i\frac{\p f^j}{\p x^\alpha}\frac{\p f^k}{\p x^\beta}\right)\ts e_i=0.\label{realharmonic}\eeq
On the other hand, $df$ is also a section of the vector bundle
$F:=T^*M\ts
  f^*(TN)$. Let $\nabla^F$ be the induced connection on $F$ by the
  Levi-Civita connections of $M$ and $N$. $f$ is said to be \emph{totally
  geodesic} if
 $ \nabla^F df=0 \in \Gamma(M,T^*M\ts T^*M\ts f^*(TN)) $, i.e.
 \beq \left(\frac{\p^2 f^i}{\p x^\alpha\p x^\beta}-\frac{\p f^i}{\p x^\gamma}\Gamma_{\alpha\beta}^\gamma
 +\Gamma_{jk}^i\frac{\p f^j}{\p x^\alpha}\frac{\p f^k}{\p x^\beta}\right)dx^\alpha\ts dx^\beta \ts e_i=0. \eeq
 Let $f:(M,h)\>(N,g)$ be an immersion, then $f$ is said to be
\emph{minimal} if $Tr_h\nabla^Fdf=0$. It is obvious that an
immersion is minimal if and only if it is harmonic.

\subsection{Harmonic maps from Hermitian manifolds to K\"ahler manifolds}

 Let $(M, h)$ be a compact Hermitian manifold and $(N,g)$ a compact
K\"ahler manifold. Let $\{z^\alpha\}_{\alpha=1}^m$ and
$\{w^i\}_{i=1}^n$ be the local holomorphic coordinates on $M$ and
$N$ respectively, where $m=\dim_\C M$ and $n=\dim_\C N$. If
$f:M\rightarrow N$ is a smooth map, the pullback
 vector bundle $f^*(T^{1,0}N)$ is denoted by $E$. The local
frames of $E$ are denoted by $e_i=f^*(\frac{\p}{\p w^i})$,
$i=1,\cdots ,n$. The metric connection on $E$ induced by the
complexified Levi-Civita connection (i.e., Chern connection) of
$T^{1,0}M$ is denoted by $\nabla^E$. There are three $E$-valued
$1$-forms, namely, \beq \bp f =\frac{\p f^i}{\p \bar z^\beta}d\bar
z^\beta\ts e_i, \ \ \p f=\frac{\p f^i}{\p z^\alpha}dz^\alpha\ts
e_i,\ \ \  df=\bp f+\p f. \eeq The $\bp$-energy of a smooth map
$f:(M,h)\>(N,g)$ is defined by \beq E''(f)=\int_M |\bp f|^2
\frac{\omega_h^m}{m!} \eeq and similarly we can define the
$\p$-energy  $E'(f)$ and the total energy $E(f)$ by \beq
E'(f)=\int_M |\p f|^2 \frac{\omega_h^m}{m!},\ \ \ \ \ E(f)=\int_M
|df|^2\frac{\omega_h^m}{m!}. \eeq It is obvious that the quantity
$E(f)$ coincides with the energy defined by the background
Riemannian metrics. The following result is well-known and  we
include a proof here  for the sake of completeness.

\blemma\label{EL3} The Euler-Lagrange equation of $\bp$-energy
$E''(f)$ is $ \bp_E^*\bp f=0 $; and the Euler-Lagrange equations of
$E'(f)$ and $E(f)$ are $ \p_E^*\p f=0$ and  $\bp^*_E\bp f+\p_E^*\p
f=0$ respectively.

 \bproof Let $F:M\times \C\>N$ be a smooth function such that
$$F(z,0)=f(z),\ \ \ \frac{\p F}{\p
t}\big|_{t=0}=v,\ \ \ \frac{\p F}{\p\bar t}\big|_{t=0}=\mu.$$ Now we
set $K= F^*\left(T^{1,0}N\right)$. The connection on $K$ induced by
the Chern connection on $T^{1,0}N$ is denoted by $\nabla^K$. Its
$(1,0)$ and $(0,1)$ components are denoted by $\p_K$ and $\bp_K$
respectively. The induced bases $F^*(\frac{\p}{\p w^i})$ of $K$ are
denoted by $\hat e_i$, $i=1,\cdots,n$.
 Since the connection
$\nabla^K$ is compatible with the Hermitian metric on $K$, we obtain
\beq \frac{\p}{\p t} E''(f_t)=\int_M\left\la \left(\p_{K}\bp
f_t\right)\left(\frac{\p}{\p t}\right), \bp f_t\right\ra
\frac{\omega_h^{m}}{m!}+\int_M \left\la \bp f_t,\left(\bp_{K}\bp
f_t\right)\left(\frac{\p}{\p \bar t}\right)\right\ra
\frac{\omega_h^{m}}{m!}.\eeq On the other hand, \beq \p_K\bp
f_t=\p_K\left(\bp f^i_t \ts \hat e_i\right)=\p_t\left(\bp
f^i_t\right)\ts \hat e_i-\bp f^i_t\wedge \nabla^{'K}\hat e_i \eeq
where $\p_t$ is $\p$-operator on the manifold $M\times \C$. By
definition, \be \nabla^{K}\hat e_i&=&F^*\left(\nabla \frac{\p}{\p
w^i}\right)=F^*\left(\Gamma_{ji}^k dz^j\ts \frac{\p}{\p
w^k}\right)=\Gamma_{ji}^kdF^j\ts \hat e_k. \ee Therefore,$\left(
\nabla^{'K}\hat e_i\right)\left(\frac{\p}{\p t}\right)=\Gamma_{ji}^k
\frac{\p F^j}{\p t}\ts \hat e_k $ and \beq \left(\p_{K}\bp
f_t\right)\left(\frac{\p}{\p t}\right)=\left(\bp\left(\frac{\p
F^i}{\p t}\right)+\bp F^k \frac{\p F^j}{\p t}\Gamma_{jk}^i\right)
\ts\hat e_i. \eeq When $t=0$, \beq \left(\p_{K}\bp
f_t\right)\left(\frac{\p}{\p t}\right)\big|_{t=0}=\left(\bp v^i+\bp
f^k\cdot v^j\cdot\Gamma_{jk}^i\right)\ts  e_i=\bp_E v\eeq since
$(N,g)$ is K\"ahler, i.e. $\Gamma_{jk}^i=\Gamma_{kj}^i$. Similarly,
we get
 \beq \left(\bp_{K}\bp
f_t\right)\left(\frac{\p}{\p \bar t}\right)\big|_{t=0}=\bp_E
\mu.\eeq Finally, we obtain \beq \frac{\p}{\p t}
E''(f_t)\big|_{t=0}=\int_M\la \bp_E v, \bp f\ra
\frac{\omega_h^m}{m!}+\int_M \la \bp f,\bp_E \mu\ra
\frac{\omega_h^m}{m!}.\eeq Hence the Euler-Lagrange equation of
$E''(f)$ is $\bp_E^*\bp f=0$. Similarly, we can get the
Euler-Lagrange equations of $E'(f)$ and $E(f)$.\eproof \elemma

For any smooth function $\Phi$ on the compact Hermitian manifold
$(M,h)$, we know \beq \begin{cases}\Delta_{\bp} \Phi=\bp^*\bp
\Phi=-h^{\alpha\bar\beta}\left(\frac{\p^2 \Phi}{\p z^\alpha\p\bar
z^\beta}-2\Gamma_{\alpha\bar\beta}^{\bar \gamma}\frac{\p
\Phi}{\p\bar z^\gamma}\right),\\ \Delta_{\p} \Phi=\p^*\p
\Phi=-h^{\alpha\bar\beta}\left(\frac{\p^2 \Phi}{\p z^\alpha\p\bar
z^\beta}-2\Gamma_{\bar\beta \alpha}^{ \gamma}\frac{\p \Phi}{\p
z^\gamma}\right),\\
\Delta_d \Phi=d^*d\Phi=\Delta_{\bp}\Phi+\Delta_\p\Phi\end{cases}
\label{b}\eeq where \beq
\Gamma_{\alpha\bar\beta}^{\gamma}=\frac{1}{2}
h^{\gamma\bar\delta}\left(\frac{\p h_{\alpha \bar\delta}}{\p\bar
z^\beta}-\frac{\p h_{\alpha \bar\beta}}{\p\bar
z^\delta}\right)=\Gamma_{\bar\beta\alpha}^\gamma=\bar{\Gamma_{\beta\bar\alpha}^{\bar\gamma}},\
\ \ \ \ \ \  \Gamma_{\alpha\beta}^\gamma=\frac{1}{2}
h^{\gamma\bar\delta}\left(\frac{\p h_{\alpha \bar\delta}}{\p
z^\beta}+\frac{\p h_{\beta \bar\delta}}{\p z^\alpha}\right).\eeq

\noindent Therefore, by Lemma \ref{computationlemma}
\begin{eqnarray}\bp_E^*\bp f\nonumber&=&\left(\bp^*\bp
f^i-h^{\alpha\bar\beta}\Gamma_{jk}^i\frac{\p f^j}{\p \bar
z^\beta}\frac{\p f^k}{\p z^\alpha}\right)\ts
e_i\\&=&-h^{\alpha\bar\beta}\left(\frac{\p^2 f^i}{\p z^\alpha\p\bar
z^\beta}-2\Gamma_{\alpha\bar\beta}^{\bar \gamma}\frac{\p f^i}{\p\bar
z^\gamma}+\Gamma_{jk}^i\frac{\p f^j}{\p \bar z^\beta}\frac{\p
f^k}{\p z^\alpha}\right)\ts e_i\label{d}\end{eqnarray}  and
\begin{eqnarray}\p_E^*\p f\nonumber&=&\left(\p^*\p
f^i-h^{\alpha\bar\beta}\Gamma_{jk}^i\frac{\p f^j}{\p \bar
z^\beta}\frac{\p f^k}{\p z^\alpha}\right)\ts e_i\\&=&
-h^{\alpha\bar\beta}\left(\frac{\p^2 f^i}{\p z^\alpha\p\bar
z^\beta}-2\Gamma_{\bar\beta\alpha}^{ \gamma}\frac{\p f^i}{\p
z^\gamma}+\Gamma_{jk}^i\frac{\p f^j}{\p \bar z^\beta}\frac{\p
f^k}{\p z^\alpha}\right)\ts e_i. \label{a}\end{eqnarray} For more
details about the computations, see e.g. \cite{LY12}.

We clarify and summarize the definitions of various harmonic maps in
the following:

\bdefinition Let $(M, h)$ be a compact Hermitian manifold and
$(N,g)$ a  K\"ahler manifold. Let  $f:(M,h)\>(N,g)$ be a smooth map
and $E=f^*(T^{1,0}N)$.

\begin{enumerate}
\item $f$ is called \emph{$\bp$-harmonic} if it is a critical point of
$\bp$-energy, i.e., $\bp^*_E\bp f=0$;

\item $f$ is called \emph{$\p$-harmonic} if it is a critical point of
$\p$-energy, i.e., $\p_E^*\p f=0$;

\item $f$ is called \emph{harmonic} if it is a critical point of
$d$-energy, i.e., $\bp_E^*\bp+\p_E^*\p f=0$, i.e. \beq
h^{\alpha\bar\beta}\left(\frac{\p^2 f^i}{\p z^\alpha\p\bar
z^\beta}-\Gamma_{\alpha\bar\beta}^{\bar \gamma}\frac{\p f^i}{\p\bar
z^\gamma}-\Gamma_{\bar\beta\alpha}^{ \gamma}\frac{\p f^i}{\p
z^\gamma}+\Gamma_{jk}^i\frac{\p f^j}{\p \bar z^\beta}\frac{\p
f^k}{\p z^\alpha}\right)\ts e_i=0;\label{harmonic}\eeq

\item $f$ is called \emph{Hermitian harmonic} if it satisfies
\beq -h^{\alpha\bar\beta}\left(\frac{\p^2 f^i}{\p z^\alpha\p\bar
z^\beta}+\Gamma_{jk}^i\frac{\p f^j}{\p \bar z^\beta}\frac{\p f^k}{\p
z^\alpha}\right)\ts e_i=0;\label{pseudo} \eeq

\item $f$ is called \emph{pluri-harmonic} if it satisfies $\p_E\bp
f=0$, i.e. \beq \frac{\p^2 f^i}{\p z^\alpha\p\bar
z^\beta}+\Gamma_{jk}^i\frac{\p f^j}{\p \bar z^\beta}\frac{\p f^k}{\p
z^\alpha}=0\label{pluriharmonic}\eeq for any $\alpha,\beta$ and $i$.

\end{enumerate}

 \edefinition

\bremark \begin{enumerate} 

\item The Hermitian harmonic equation (\ref{pseudo}) was firstly
introduced by Jost-Yau in \cite{JY}. For more generalizations, see
\cite{LE}, \cite{KG} and also the references therein;

\item The harmonic map equation (\ref{harmonic}) is the same as
 classical harmonic equation (\ref{realharmonic}) by using the background Riemmanian metrics;

\item Pluri-harmonic maps are Hermitian harmonic;

\item Pluri-harmonic maps are not necessarily $\p$-harmonic or $\bp$-harmonic;

\item $\bp_E\p f=0$ and $\p_E\bp f=0$ are equivalent;

\item For another type of Hermitian harmonic maps between Hermitian
manifolds defined by using Chern connections, we refer the reader to
\cite{ZX}.
\end{enumerate}

\eremark

\blemma For any smooth map $f:(M,h)\>(N,g)$ from a Hermitian
manifold $(M,h)$ to a K\"ahler manifold $(N,g)$, we have $\bp_E\bp
f=0$ and $\p_E\p f=0. $ \bproof It is easy to see that \be \bp_E\bp
f&=&\bp_E\left(\frac{\p f^i}{\p
\bar z^\alpha}d\bar z^\alpha\ts e_i\right)\\
&=&\frac{\p^2 f^i}{\p\bar z^\beta\p\bar z^\alpha}d\bar z^\beta\wedge
d\bar z^\alpha \ts e_i-\frac{\p f^i}{\p\bar z^\alpha}d\bar
z^\alpha\wedge \Gamma_{ik}^j\frac{\p f^k}{\p\bar z^\beta}d\bar
z^\beta\ts e_j=0\ee since $\Gamma_{ik}^j=\Gamma_{ki}^j$ when $(N,h)$
is K\"ahler. The proof of $\p_E\p f=0$ is similar. \eproof \elemma

\blemma\label{pseudo1} Let $ f:(M,h)\>(N,g)$ be a smooth map from a
compact Hermitian manifold $(M,h)$ to a compact K\"ahler manifold
$(N,g)$. The Hermitian harmonic map equation (\ref{pseudo}) is
equivalent to
 \beq
\left(\bp_E-2\sq\p^*\omega_h\right)^*\left(\bp f\right)=0 \qtq{or}
\left(\p_E+2\sq\bp^*\omega_h\right)^*\left(\p f\right)=0. \eeq
\bproof On a compact Hermitian manifold $(M,h)$ with
$\omega_h=\frac{\sq}{2}h_{\alpha\bar\beta}dz^\alpha\wedge d\bar
z^\beta$, we have(\cite[Lemma~A.6]{LY12}) \beq \bp^*\omega_h=\sq
\Gamma_{\gamma \bar \beta}^{\bar\beta}dz^\gamma \qtq{and} -2\sq
(\bp^*\omega_h)^*(\bp
f)=-2h^{\alpha\bar\beta}\Gamma_{\alpha\bar\beta}^{\gamma}\frac{\p
f^i}{\p z^\gamma}\ts e_i.\label{c}\eeq The equivalence is derived
from (\ref{c}),  (\ref{a}) and (\ref{d}).
 \eproof
\elemma

\bdefinition A compact Hermitian manifold $(M,h)$ is call
\emph{balanced} if the fundamental form $\omega_h$ is co-closed,
i.e. $d^*\omega_h=0$. \edefinition

\bproposition\label{balanced} Let $(M,h)$ be a compact balanced
Hermitian manifold and $(N,g)$ a  K\"ahler manifold. The $E'$, $E''$
and $E$-critical points coincide. Moreover, they satisfy the
Hermitian harmonic equation (\ref{pseudo}). That is, $\bp$-harmonic,
$\p$-harmonic, Hermitian harmonic and harmonic maps are the same if
the domain $(M,h)$ is a balanced manifold.

 \bproof The
balanced condition $d^*\omega_h=0$ is equivalent to $\p^*\omega_h=0$
or $\bp^*\omega_h=0$ or $
h^{\alpha\bar\beta}\Gamma_{\alpha\bar\beta}^{\gamma}=0 $ for
$\gamma=1,\cdots, m.$ By formulas (\ref{a}) and (\ref{d}), we obtain
$\bp_E^*\bp f=\p_E^*\p f$. The second statement follows by Lemma
\ref{pseudo1}. \eproof \eproposition

\bproposition Let $(M,h)$ be a compact Hermitian manifold and
$(N,g)$ a  K\"ahler manifold.  If $f:(M,h)\>(N,g)$ is
  totally geodesic and $(M,h)$ is K\"ahler, then $f$ is
  pluri-harmonic.

 \bproof Considering the complexified connection $\nabla^F$ on the vector bundle $F=T^*M\ts
  f^*(TN)$, we have \be
\nabla^F \bp f&=&\nabla^F\left(\frac{\p f^i}{\p \bar z^\beta}d\bar
z^\beta\ts
e_i\right)\\
&=&\left(\frac{\p^2 f^i}{\p z^\alpha\p\bar z^\beta}-\frac{\p
f^i}{\p\bar z^\gamma}\Gamma_{\alpha\bar\beta}^{\bar \gamma}+\frac{\p
f^j}{\p z^\alpha}\frac{\p f^k}{\p \bar z^\beta}\Gamma_{jk}^i\right)
dz^\alpha\ts d\bar z^\beta\ts e_i\\
&&+\left(\frac{\p^2 f^i}{\p\bar z^\gamma\p\bar z^\delta}-\frac{\p
f^i}{\p\bar z^\alpha}\Gamma_{\bar\gamma\bar\delta}^{\bar
\alpha}+\frac{\p f^j}{\p\bar z^\gamma}\frac{\p f^k}{\p\bar
z^\delta}\Gamma_{jk}^i\right)d\bar z^\gamma\ts d\bar z^\delta\ts
e_i\\&&-\frac{\p f^i}{\p\bar
z^\alpha}\Gamma_{\bar\beta\lambda}^{\bar\alpha}d\bar z^\beta\ts d
z^\lambda \ts e_i. \ee

\noindent Similarly, we have \be \nabla^F \p
f&=&\nabla^F\left(\frac{\p f^i}{\p  z^\alpha}d z^\alpha\ts
e_i\right)\\
&=&\left(\frac{\p^2 f^i}{\p z^\alpha\p\bar z^\beta}-\frac{\p f^i}{\p
z^\gamma}\Gamma_{\bar\beta\alpha}^{ \gamma}+\frac{\p f^j}{\p
z^\alpha}\frac{\p f^k}{\p \bar z^\beta}\Gamma_{jk}^i\right)
 d\bar z^\beta\ts dz^\alpha\ts e_i\\
&&+\left(\frac{\p^2 f^i}{\p z^\gamma\p z^\delta}-\frac{\p f^i}{\p
z^\alpha}\Gamma_{\gamma\delta}^{ \alpha}+\frac{\p f^j}{\p
z^\gamma}\frac{\p f^k}{\p z^\delta}\Gamma_{jk}^i\right)d z^\gamma\ts
d z^\delta\ts e_i\\&&-\frac{\p f^i}{\p
z^\alpha}\Gamma_{\lambda\bar\beta}^{\alpha}d z^\lambda \ts d\bar
z^\beta\ts e_i. \ee

\noindent
 That is \be\nabla^F df &=&\nabla^F\bp f+\nabla^F\p f\\&=&\left(\frac{\p^2 f^i}{\p
z^\alpha\p\bar z^\beta}-\frac{\p f^i}{\p\bar
z^\gamma}\Gamma_{\alpha\bar\beta}^{\bar \gamma}-\frac{\p f^i}{\p
z^\lambda}\Gamma_{\alpha\bar\beta}^\lambda+\frac{\p f^j}{\p
z^\alpha}\frac{\p f^k}{\p \bar z^\beta}\Gamma_{jk}^i\right)
dz^\alpha\ts d\bar z^\beta\ts e_i\\
&&+\left(\frac{\p^2 f^i}{\p\bar z^\gamma\p\bar z^\delta}-\frac{\p
f^i}{\p\bar z^\alpha}\Gamma_{\bar\gamma\bar\delta}^{\bar
\alpha}+\frac{\p f^j}{\p\bar z^\gamma}\frac{\p f^k}{\p\bar
z^\delta}\Gamma_{jk}^i\right)d\bar z^\gamma\ts d\bar z^\delta\ts
e_i\\&&+\left(\frac{\p^2 f^i}{\p z^\alpha\p\bar z^\beta}-\frac{\p
f^i}{\p z^\gamma}\Gamma_{\bar\beta\alpha}^{ \gamma}-\frac{\p
f^i}{\p\bar z^\delta}\Gamma_{\bar\beta\alpha}^{\bar \delta}+\frac{\p
f^j}{\p z^\alpha}\frac{\p f^k}{\p \bar z^\beta}\Gamma_{jk}^i\right)
 d\bar z^\beta\ts dz^\alpha\ts e_i\\ &&+\left(\frac{\p^2 f^i}{\p z^\gamma\p z^\delta}-\frac{\p
f^i}{\p z^\alpha}\Gamma_{\gamma\delta}^{ \alpha}+\frac{\p f^j}{\p
z^\gamma}\frac{\p f^k}{\p z^\delta}\Gamma_{jk}^i\right)d z^\gamma\ts
d z^\delta\ts e_i \ee

\noindent If $f$ is totally geodesic and $(M,h)$ is K\"ahler, then
$f$ is pluri-harmonic by degree reasons. \eproof \eproposition

\bremark It is easy to see that pluri-harmonic maps are not
necessarily totally geodesic. \eremark

\blemma\label{key111} Let $f$ be a pluri-harmonic map from a complex
manifold $M$ to a K\"ahler manifold $(N,g)$. Then the real $(1,1)$
forms \beq \omega_0=\frac{\sq}{2}g_{i\bar j}\frac{\p f^i}{\p
z^\alpha}\frac{\p \bar f^j}{\p \bar z^\beta} dz^\alpha\wedge d\bar
z^\beta=\frac{\sq}{2}g_{i\bar j} \p f^i\wedge \bp \bar f^j ,\eeq and
\beq \omega_1=\frac{\sq}{2}g_{i\bar j}\frac{\p \bar f^j}{\p
z^\alpha}\frac{\p f^i}{\p \bar z^\beta}dz^\alpha\wedge d\bar
z^\beta= \frac{\sq}{2} g_{i\bar j}\p\bar f^i\wedge \bp f^j\eeq
 are all
$d$-closed, i.e.,
$$d\omega_0=\p\omega_0=\bp\omega_0=0,\qtq{and}d\omega_1=\p\omega_1=\bp\omega_1=0.$$
\bproof By definition, we see \be \p\omega_0&=&-\frac{\sq}{2}\p
f^i\wedge \p\left(g_{i\bar j} \bp \bar{f}^{j}\right)\\
&=&-\frac{\sq}{2}\p f^i\wedge \left(\frac{\p g_{i\bar j}}{\p
z^k}\cdot\p f^k\wedge \bp\bar f^j+\frac{\p g_{i\bar j}}{\p \bar
z^\ell}\cdot\p \bar f^\ell\wedge \bp\bar f^j+g_{i\bar j} \p\bp
\bar{f}^{j}\right)\\
(\text{$g$ is K\"ahler})&=&-\frac{\sq}{2}\p f^i\wedge \left(\frac{\p
g_{i\bar j}}{\p \bar z^\ell}\cdot\p \bar f^\ell\wedge \bp\bar
f^j+g_{i\bar j} \p\bp \bar{f}^{j}\right)\\
&=&-\frac{\sq}{2}\p f^i\wedge  g_{i\bar
s}\left(\p\bp\bar{f}^s+g^{p\bar s}\frac{\p g_{p\bar q}}{\p \bar
z^\ell}\cdot\p \bar f^\ell\wedge \bp\bar f^q\right)\\
&=&-\frac{\sq}{2}\p f^i\wedge  g_{i\bar
s}\left(\p\bp\bar{f}^s+\bar{\Gamma_{q\ell}^s}\cdot\p \bar
f^\ell\wedge \bp\bar f^q\right)\\
&=&0\ee where the last step follows from the definition equation
(\ref{pluriharmonic}) of pluri-harmonic maps. Hence, we obtain
$d\omega_0=0$. The proof of $d\omega_1=0$ is similar.
 \eproof
\elemma

\subsection{Harmonic maps between  Hermitian manifolds}

Let $(M,h)$ and $(N,g)$ be two compact Hermitian manifolds. Using
the same notation as in the previous subsection, we can define
$\bp$-harmonic (resp. $\p$-harmonic, harmonic ) map $f:(M,h)\>(N,g)$
by using the critical point of the Euler-Lagrange equation of
$E''(f)$ (resp. $E'(f)$, $E(f)$). In this case, the harmonic
equations have the same second order parts, but the torsion parts
are different. For example, the $\bp$-harmonic equation is \beq
\left(\Delta_{\bp} f^i+ T^i(f)\right)\ts e_i=0 \eeq where $T(f)$ is
a quadratic function in $df$ and the coefficients are the
Christoffel symbols  of $(N,g)$. One can see it  clearly from the
proof of Lemma \ref{EL3}.

\subsection{Harmonic maps from Hermitian manifolds to Riemannian manifolds}
Let $(M,h)$ be a compact Hermitian manifold, $(N,g)$ a Riemannian
manifold and $E=f^*(TN)$ with the induced Levi-Civita connection. As
similar as in the K\"ahler target manifold case, we can define the
$\bp$-energy of $f:(M,h)\>(N,g)$ \beq E''(f)=\int_M|\bp f|^2
\frac{\omega_h^m}{m!}=\int_{M} g_{ij}h^{\alpha\bar\beta}\frac{\p
f^i}{\p z^\alpha}\frac{\p f^j}{\p\bar
z^\beta}\frac{\omega_h^m}{m!}.\label{realpartial}\eeq It is easy to
see that the Euler-Lagrange equation of (\ref{realpartial}) is
\beq\bp_E^*\bp f= \Delta_{\bp}
f^i-h^{\alpha\bar\beta}\Gamma_{jk}^i\frac{\p f^j}{\p
z^\alpha}\frac{\p f^k}{\p \bar z^\beta}=0. \eeq Similarly, we can
define  $E'(f)$ and get its Euler-Lagrange equation \beq \p_E^*\p
f=\Delta_{\p} f^i-h^{\alpha\bar\beta}\Gamma_{jk}^i\frac{\p f^j}{\p
z^\alpha}\frac{\p f^k}{\p \bar z^\beta}=0. \eeq The Euler-Lagrange
equation of $E(f)$ is $ \bp_E^*\bp f+\p_E^*\p f=0$.

\bdefinition Let $(M, h)$ be a compact Hermitian manifold and
$(N,g)$ a  Riemannian manifold. Let  $f:(M,h)\>(N,g)$ be a smooth
map and $E=f^*(TN)$.
\begin{enumerate}
\item $f$ is called \emph{$\bp$-harmonic} if it is a critical point of
$\bp$-energy, i.e., $\bp^*_E\bp f=0$;

\item $f$ is called \emph{$\p$-harmonic} if it is a critical point of
$\p$-energy, i.e., $\p_E^*\p f=0$;

\item $f$ is called \emph{harmonic} if it is a critical point of
$d$-energy, i.e., $\bp_E^*\bp+\p_E^*\p f=0$;

\item $f$ is called \emph{Hermitian harmonic} if it satisfies
\beq -h^{\alpha\bar\beta}\left(\frac{\p^2 f^i}{\p z^\alpha\p\bar
z^\beta}+\Gamma_{jk}^i\frac{\p f^j}{\p \bar z^\beta}\frac{\p f^k}{\p
z^\alpha}\right)\ts e_i=0;\label{pseudo2} \eeq

\item $f$ is called \emph{pluri-harmonic} if it satisfies $\p_E\bp
f=0$, i.e. \beq \left(\frac{\p^2 f^i}{\p z^\alpha\p\bar
z^\beta}+\Gamma_{jk}^i\frac{\p f^j}{\p \bar z^\beta}\frac{\p f^k}{\p
z^\alpha}\right)dz^\alpha\wedge d\bar z^\beta\ts e_i=0.\eeq
\end{enumerate} \edefinition

\noindent As similar as Proposition \ref{balanced}, we have

\bcorollary Let $f:(M,h)\>(N,g)$ be a smooth map from a compact
Hermitian manifold $(M,h)$ to a  Riemannian manifold $(N,g)$. If
$(M,h)$ is a balanced Hermitian manifold, i.e. $d^*\omega_h=0$, then
$\p$-harmonic map, $\bp$-harmonic map, Hermitian harmonic map and
harmonic map coincide. \ecorollary


\section{Manifolds with non-degenerate curvatures}

\subsection{Curvatures of K\"ahler manifolds}

Let $(N,g)$ be a K\"ahler manifold. In the local holomorphic
coordinates $(w^1,\cdots, w^n)$ of $N$, the curvature tensor
components are \beq R_{i\bar j k\bar \ell}=-\frac{\p^2 g_{k\bar
\ell}}{\p w^i\p\bar w^j}+g^{p\bar q}\frac{\p g_{k\bar q}}{\p
w^i}\frac{\p g_{p\bar \ell}}{\p\bar w^j}. \eeq In \cite{Siu}, Siu
introduced the following definition:  the curvature tensor $R_{i\bar
j k\bar \ell}$ is said to be \emph{strongly negative} (resp.
\emph{strongly positive}) if \beq \sum_{i,j,k,\ell}R_{i\bar j
k\bar\ell}
(A^{i}\overline{B}^{j}-C^{i}\overline{D}^{j})(\overline{A^{\ell}\overline{B}^{k}-C^{\ell}\overline{D}^{k}})<0,
\qtq{(resp. $>$0)}\label{strong}\eeq for any nonzero $n\times n$
complex matrix
$(A^{i}\overline{B}^{j}-C^{i}\overline{D}^{j})_{i,j}$.

\bdefinition\label{Siudef} Let $(N,g)$ be a  K\"ahler manifold. The
curvature tensor $R_{i\bar j k\bar \ell}$ is called
\emph{non-degenerate} if it satisfies the condition that \beq
\sum_{i,j,k,\ell}R_{i\bar j k\bar\ell}
(A^{i}\overline{B}^{j}-C^{i}\overline{D}^{j})(\overline{A^{\ell}\overline{B}^{k}-C^{\ell}\overline{D}^{k}})=0\label{nondegenerate}\eeq
if and only if $A^{i}\overline{B}^{j}-C^{i}\overline{D}^{j}=0$ for
any $i,j$. \edefinition \noindent It is easy to see that both
manifolds with strongly positive curvatures and manifolds with
strongly negative curvatures have non-degenerate curvatures.

\subsection{Curvatures of Riemannian manifolds}
 Let $(M,g)$ be a Riemannian manifold and $\nabla$ the Levi-Civita
  connection. The curvature tensor $R$ is defined by
  $$R(X,Y,Z,W)=g(\nabla_X\nabla_YZ-\nabla_Y\nabla_XZ-\nabla_{[X,Y]}Z, W).$$
To illustrate our computation rules on Riemannain manifolds, for
example, the Riemannian curvature tensor components of $S^n$ induced
by the canonical metric of $\R^{n+1}$ are $
R_{ijk\ell}=g_{i\ell}g_{jk}-g_{ik}g_{j\ell}. $ The Ricci curvature
tensor components are $ R_{jk}=g^{i\ell}R_{ijk\ell}=(n-1)g_{jk}.
$

As similar as Siu's definition, Sampson(\cite{Sa2}) proposed the
following definition:

\bdefinition \label{samdef}Let $(M,g)$ be a compact Riemannian
manifold. \bd\item The curvature tensor $R$ of $(M,g)$ is said to be
\emph{Hermitian-positive} (resp. \emph{Hermitian-negative}) if \beq
R_{ijk\ell}A^{i\bar \ell}A^{j\bar k}\geq 0 (\qtq{resp.} \leq 0)
\label{hermitiannegative}\eeq for any Hermitian semi-positive matrix
$A=(A^{i\bar \ell})$. $R$ is called \emph{strongly
Hermitian-positive} (resp. \emph{strongly Hermitian-negative}) if
$R$ is Hermitian-positive (resp. Hermitian negative) and the
equality in (\ref{hermitiannegative}) holds only for Hermitian
semi-positive matrix $A$ with complex rank $\leq 1$.

\noindent
\item $R$ is said to be \emph{Hermitian non-degenerate at some point $p\in M$} if  \beq R_{ijk\ell}(p)A^{i\bar \ell}A^{j\bar
k}=0\eeq for some Hermitian semi-positive matrix $A=(A^{i\bar j})$
implies $A$ has rank $\leq 1$. $R$ is said to be \emph{Hermitian
non-degenerate} if it is Hermitian non-degenerate everywhere. \ed
\edefinition

Note that any rank one Hermitian matrix  can be written as $A^{i\bar
j}=a^ib^{\bar j}$ and so for any curvature tensor $R_{ijk\ell}$, one
has
$$R_{ijk\ell}A^{i\bar
\ell}A^{j\bar k}=0.$$ On the other hand, it is easy to see that both
manifolds with strongly Hermitian-positive curvatures and manifolds
with strongly Hermitian-negative curvatures are  Hermitian
non-degenerate . \blemma[\cite{Sa2}] If $(M,g)$ has positive (resp.
negative) constant sectional curvature, then the curvature tensor is
strongly Hermitian positive (resp. negative).  In particular, it is
Hermitian non-degenerate. \bproof Let
$R_{ijk\ell}=\kappa(g_{i\ell}g_{jk}-g_{ik}g_{j\ell})$. Then \beq
R_{ijk\ell}A^{i\bar \ell}A^{j\bar
k}=\kappa\left((TrA)^2-Tr(A^2)\right). \eeq The results follow by
this identity easily. \eproof \elemma

\bremark\label{WP} In \cite{LSYY}, we give a complete list on the
curvature relations  of a K\"ahler manifold $(M,g)$: \bd
\item semi dual-Nakano-negative;

\item non-positive Riemannian curvature operator;

\item strongly non-positive in the sense of siu;

\item non-positive complex sectional curvature;

\item non-positive Riemannian sectional curvature;

\item non-positive holomorphic bisectional curvature;

\item non-positive isotropic curvature.
$$ (1)\Longrightarrow(2)\Longrightarrow (3)\Longleftrightarrow(4)\Longrightarrow (5)\Longrightarrow (6)$$
$$(1)\Longrightarrow (3)\Longleftrightarrow(4)\Longrightarrow (7).$$\ed
\noindent So far, it is not clear to the authors whether one of them
can imply (Sampson's) Hermitian negativity. However, it is easy to
see that the  Poincar\'e disks and projective spaces have  Hermitian
negative and Hermitian positive curvatures respectively. It is
hopeful that semi dual-Nakano-negative curvatures can imply
Hermitian negative curvatures (in the sense of Sampson). We will go
back to this topic later. \eremark

\section{Existence  of various harmonic maps}\label{existence}

In their pioneering work \cite{ES}, Eells-Sampson have proposed the
heat flow method to study the existence of harmonic maps. In this
section, we will consider a similar setting. Let $f:(M,h)\>(N,g)$ be
a continuous map from a compact Hermitian manifold to a compact
Riemannian manifold. In the paper (\cite{JY}) of Jost and Yau, they
considered the heat flow for the Hermitian harmonic equation, i.e.,
\beq
\begin{cases}
h^{\alpha\bar\beta}\left(\frac{\p^2 f^i(z,t)}{\p z^\alpha\p\bar
z^\beta}+\Gamma_{jk}^i\frac{\p f^j(z,t)}{\p \bar z^\beta}\frac{\p
f^k(z,t)}{\p z^\alpha}\right)-\frac{\p f^i(z,t)}{\p
t}=0\\
f_0=f
\end{cases}\label{JY1}
\eeq where $\Gamma_{jk}^i$ are Christoffel symbols of the Riemannian
manifold $(N,g)$.\blemma[Jost-Yau] If $(N,g)$ has non-positive
Riemannian sectional curvature, then a solution of (\ref{JY1})
exists for all $t\geq 0$.

\elemma \noindent Similarly,
 we can consider the following parabolic system
for  the $\bp$-energy of a smooth map $f$ from a compact Hermitian
manifolds $(M,h)$ to a Riemannian manifold $(N,g)$, \beq
\begin{cases}
\frac{d f_t}{dt}=-\bp_E^*\bp f_t\\
f_0=f.
\end{cases}
\label{heatflow}\eeq

\noindent Locally, the parabolic equation (\ref{heatflow}) is \beq
h^{\alpha\bar\beta}\left(\frac{\p^2 f^i(z,t)}{\p z^\alpha\p\bar
z^\beta}-2\Gamma_{\alpha\bar\beta}^{\bar \gamma}\frac{\p
f^i(z,t)}{\p\bar z^\gamma}+\Gamma_{jk}^i\frac{\p f^j(z,t)}{\p \bar
z^\beta}\frac{\p f^k(z,t)}{\p z^\alpha}\right)-\frac{\p f^i(z,t)}{\p
t}=0.\label{P}\eeq
The difference between (\ref{P}) and (\ref{JY1}) are the first order
derivative terms of $f$. By the theory of parabolic PDEs, if $(N,g)$
has non-positive sectional curvature, the solution of (\ref{P})
exists for all $t\geq 0$ following the adapted methods in \cite{ES}
and \cite{JY}. Let
$$e(f)=h^{\alpha\bar\beta}g_{ij}\frac{\p f^i}{\p z^\alpha}\frac{\p
f^j}{\p\bar z^\beta}$$ be the energy density. By differentiating the
equation (\ref{P}), we obtain
  \beq
\left(\Delta_c-\frac{\p}{\p t}\right)e(f)\geq \frac{1}{2}|\nabla^2
f|^2 -C e(f)\eeq if $(N,g)$ has non-positive sectional curvature
where $C=C(M,h)$ is a positive constant only depending  on $(M,h)$,
and $\Delta_c$ is the canonical Laplacian
$\Delta_c=h^{\alpha\bar\beta}\frac{\p^2}{\p z^\alpha\p\bar
z^\beta}.$ The extra first order terms in $f$ are absorbed in
$|\nabla^2 f|$ by using the Schwarz inequality. As analogous to the
existence results of Eells-Sampson(\cite{ES}, harmonic maps) and
Jost-Yau(\cite{JY}, Hermitian harmonic maps), we obtain

\btheorem\label{realexistence} Let $(M,h)$ be a compact Hermitian
manifold and $(N,g)$  a compact Riemannian manifold of negative
Riemannian sectional curvature. Let $\phi:M\>N$ be continuous, and
suppose that $\phi$ is not homotopic to a map onto a closed geodesic
of $N$. Then there exists a $\bp$-harmonic (resp. $\p$-harmonic) map
which  is homotopic to $\phi$.
 \etheorem

\btheorem\label{realexistence2} Let $(M,h)$ be a compact Hermitian
manifold $(N,g)$  a compact Riemannian manifold of negative
Riemannian sectional curvature. Let $\phi:M\>N$ be a continuous map
with $e(\phi^*(TN))\neq 0$ where $e$ is the Euler class. Then there
exists a $\bp$-harmonic (resp. $\p$-harmonic) map which  is
homotopic to $\phi$.
 \etheorem

\noindent As a special case, we have

\bcorollary\label{complexexistence} Let $(M,h)$ be a compact
Hermitian manifold and $(N,g)$ be a compact K\"ahler manifold of
strongly negative curvature. Let $\phi:M\>N$ be  a continuous map
and suppose that $\phi$ is not homotopic to a map onto a closed
geodesic of $N$. Then there exists a $\bp$-harmonic (resp.
$\p$-harmonic) map which  is homotopic to $\phi$. \bproof It follows
from the fact that if a K\"ahler manifold has strongly negative
curvature, then the background Riemannain metric has negative
sectional curvature.
 \eproof\ecorollary

\noindent Finally, we need to point out that, along the same line,
one can easily obtain similar existence results for various harmonic
maps into a Hermitian target manifold $(N,g)$ if the background
Riemannian metric on $N$ has non-positive Riemannian sectional
curvature. \noindent For more details about the existence and
uniqueness results on various harmonic maps in the Hermitian
context, we refer the reader to \cite{ES}, \cite{Hart}, \cite{EL},
\cite{EL2}, \cite{JY}, \cite{GK}, \cite{JS} \cite{Ni}, \cite{LY08},
\cite{Dong}, \cite{YHD} and also the references therein.

\section{The complex analyticity of harmonic maps }\label{hdd}

\subsection{Harmonic maps from Hermitian manifolds to K\"ahler
manifolds}\label{siudd} 

Let $f:(M,h)\>(N,g)$ be a smooth map from a Hermitian manifold
$(M,h)$ to a K\"ahler manifold. Let $E=f^*(T^{1,0}N)$. In the local
coordinates $\{z^\alpha\}$ on $M$, and $\{w^i\}$ on $N$,  one can
get
\begin{eqnarray} Q:&=&\sq \left\la[R^E,\Lambda]\bp f, \bp
f\right\ra\nonumber\\
&=&-\frac{1}{2}\sum_{\alpha,\gamma}R_{i\bar j k\bar
\ell}\left(\frac{\p f^i}{\p z^\alpha}\frac{\p \bar f^j}{\p
z^\gamma}-\frac{\p f^i}{\p z^\gamma}\frac{\p \bar f^j}{\p
z^\alpha}\right)\bar{\left(\frac{\p f^\ell}{\p z^\alpha}\frac{\p
\bar f^k}{\p z^\gamma}-\frac{\p f^\ell}{\p z^\gamma}\frac{\p \bar
f^k}{\p z^\alpha}\right)}\label{curvature}\end{eqnarray} in the
local normal coordinates
$h_{\alpha\bar\beta}=\delta_{\alpha\bar\beta}$ centered at a point
$p\in M$ where $R^E$ is the $(1,1)$ component of the curvature
tensor of $E$ and $R_{i\bar j k\bar \ell}$ are components of the
curvature tensor of $(N,g)$. If $Q$ is zero, $(N,g)$ has
\emph{non-degenerate curvature}, $N$ is compact and $rank_{\R}
df\geq 4$, one can show $\p f=0$ or $\bp f=0$(c.f. \cite[Siu]{Siu}).

Now let's recall Siu's $\p\bp$ trick(\cite{Siu,Siu2}) in the
Hermitian setting (c.f.\cite{JY}). Let $f:(M,h)\>(N,g)$ be a smooth
map between Hermitian manifolds and $E=f^*(T^{1,0}N)$.

\blemma We have the following formula \beq \p\bp\{\bp f, \bp
f\}=-\{\p_E\bp f,\p_E\bp f\}+\{\bp f, R^E \bp f\}. \label{siu1}\eeq
\elemma

 \blemma\label{new} Let $E$ be any Hermitian vector bundle over a
Hermitian manifold $(M,\omega),$ and $\phi$ a smooth $E$-valued
$(1,1)$-form on $M$. One has
\beq-\{\phi,\phi\}\frac{\omega^{m-2}}{(m-2)!}=4\left(
|\phi|^2-|Tr_\omega\phi|^2\right)\frac{\omega^m}{m!}.
\label{siu2}\eeq \bproof  Without loss of generality, we can assume
$E$ is a trivial bundle, and  $h_{i\bar j}=\delta_{i\bar j}$ at a
fixed point $p\in M$, then for $\phi=\phi_{p\bar q}dz^p\wedge d\bar
z^q$.
 \be
&&\left(\frac{\sq}{2}\right)^2\{\phi,\phi\}\frac{\omega^{m-2}}{(m-2)!}\\&=&\left(\frac{\sq}{2}\right)^2\left(\sum_{p,q,s,t}\phi_{p\bar
q}dz^p\wedge d\bar z^q \cdot \bar{\phi}_{s\bar t}d\bar z^s\wedge
dz^t\right)\frac{\omega^{m-2}}{(m-2)!}\\
&=&\left(\sum_{1\leq p<q\leq m}2|\phi_{p\bar q}|^2-2\sum_{1\leq
s<t\leq m}\phi_{s\bar s}\bar{\phi_{t\bar
t}}\right)\frac{\omega^m}{m!}\\&=&\left(\sum_{p,q}|\phi_{p\bar
q}|^2-\sum_{s,t}\phi_{s\bar s}\bar{\phi_{t\bar
t}}\right)\frac{\omega^m}{m!}\\
&=&\left(|\phi|^2-|Tr_\omega\phi|^2\right)\frac{\omega^m}{m!}. \ee

 \eproof
 \elemma

\bremark The right hand side of (\ref{siu2}) is not positive in
general. When $\phi$ is primitive, i.e. $Tr_\omega\phi=0$,
(\ref{siu2}) is the Riemann-Hodge bilinear relation for primitive
$(1,1)$ forms (e.g. \cite[Corollary ~1.2.36]{Huybrechets05} or
\cite[Proposition ~6.29]{Voisin02}). \eremark

\blemma We have the following formula for any smooth map $f$ from a
Hermitian manifold $(M,h)$ to a K\"ahler manifold $(N,g)$. \beq\{\bp
f, R^E\bp
f\}\frac{\omega^{m-2}_h}{(m-2)!}=4Q\cdot\frac{\omega_h^m}{m!}
\label{siu3}\eeq where $Q$ is defined in (\ref{curvature}). \elemma

\blemma\label{identity} Let $f$ be any \bf{smooth map} from a
compact Hermitian manifold $(M,h)$ to a K\"ahler manifold $(N,g)$.
We have the following identity\beq \int_M\p\bp\{\bp f, \bp f\}
\frac{\omega_h^{m-2}}{(m-2)!}=4\int_M \left(|\p_E\bp f|^2-|Tr_\omega
\p_E\bp
f|^2\right)\frac{\omega^m_h}{m!}+\int_M4Q\cdot\frac{\omega_h^m}{m!}.
\label{siu5}\eeq \bproof It follows by formula (\ref{siu1}),
(\ref{siu2}) and (\ref{siu3}). \eproof
 \elemma

Now one can  get the following generalization of Siu's result
(\cite{Siu}):

\bcorollary[{\cite[Jost-Yau]{JY}}] Let  $(N,g)$ be a compact
K\"ahler manifold, and $(M,h)$ a compact Hermitian manifold with
$\p\bp\omega^{m-2}_h=0$ where $m=\dim_\C M$. Let
 $f:(M,h)\>(N,g)$ be a Hermitian harmonic map. Then $f$ is holomorphic or
anti-holomorphic if $(N,g)$ has strongly negative curvature (in the
sense of Siu) and $rank_{\R} df\geq 4$.
 \ecorollary

\subsection{Harmonic maps from Hermitian manifolds to Riemannian
manifolds} In this subsection, we shall apply similar ideas in
Section \ref{siudd} to harmonic maps from Hermitian manifolds to
Riemannian manifolds. Let $f:(M,h)\>(N,g)$ be a smooth map from a
compact Hermitian manifold $(M,h)$ to a Riemannian manifold $(N,g)$.

\blemma\label{11pullback} The $(1,1)$-part of the curvature tensor
of $E=f^*(TN)$ is \beq R_{1,1}^{f^*(TN)}=2R_{ijk}^{\ell}\frac{\p
f^i}{\p z^\alpha}\frac{\p f^j}{\p \bar z^\beta}dz^\alpha\wedge d\bar
z^\beta\ts e^k\ts e_\ell. \eeq

 \bproof Since the  curvature tensor of the real vector bundle $TN$
is \beq R^{TN}=R_{ijk}^\ell dx^i\wedge  dx^j\ts \left(dx^k \ts
\frac{\p}{\p x^\ell}\right)\in \Gamma(N, \Lambda^2T^*N\ts End(TN)),
\eeq we get the full curvature tensor of the pullback vector bundle
$E=f^*(TN)$, \beq f^*\left(R^{TN}\right)=R_{ijk}^\ell df^i\wedge
df^j \ts e^k\ts e_\ell\in\Gamma(N,\Lambda^2 T^*M\ts End(E)). \eeq
The $(1,1)$ part of it is \be
R_{1,1}^{f^*(TN)}&=&R_{ijk}^{\ell}\left(\frac{\p f^i}{\p
z^\alpha}\frac{\p f^j}{\p \bar z^\beta}-\frac{\p f^i}{\p\bar
z^\beta}\frac{\p f^j}{\p z^\alpha}\right)dz^\alpha\wedge d\bar
z^\beta\ts e^k\ts e_\ell\\&=&2R_{ijk}^{\ell}\frac{\p f^i}{\p
z^\alpha}\frac{\p f^j}{\p \bar z^\beta}dz^\alpha\wedge d\bar
z^\beta\ts e^k\ts e_\ell, \ee since $R_{ijk}^\ell=-R_{jik}^\ell$.
\eproof
 \elemma

\blemma We have \beq \left\la \sq [R_{1,1}^{f^*(TN)},\Lambda]\p f,
\p f\right\ra=2 h^{\alpha\bar\delta} h^{\gamma\bar\beta}R_{ijk\ell}
\frac{\p f^i}{\p z^\alpha}\frac{\p f^k}{\p \bar z^\beta}\frac{\p
f^j}{\p z^\gamma}\frac{\p f^\ell}{\p \bar z^\delta}
\label{sampson1}\eeq and \beq  Q_0:=\left\la \sq
[R_{1,1}^{f^*(TN)},\Lambda]\bp f, \bp f\right\ra=-2
h^{\alpha\bar\delta} h^{\gamma\bar\beta}R_{ijk\ell} \frac{\p f^i}{\p
z^\alpha}\frac{\p f^k}{\p \bar z^\beta}\frac{\p f^j}{\p
z^\gamma}\frac{\p f^\ell}{\p \bar z^\delta}. \label{sampson2} \eeq
\bproof It is easy to see that the identity (\ref{sampson2}) is the
complex conjugate of  (\ref{sampson1}). By Lemma \ref{11pullback},
\be  \sq [R_{1,1}^{f^*(TN)},\Lambda]\p f&=&-\sq\Lambda R^{E}\p
f\\&=& 2h^{\alpha\bar\beta}\left(-R_{ijk}^\ell+R_{kji}^\ell\right)
\frac{\p f^i}{\p z^\alpha}\frac{\p f^j}{\p \bar z^\beta}\frac{\p
f^k}{\p
z^\gamma}dz^\gamma\ts e_\ell\\
&=& 2h^{\alpha\bar\beta}R_{kij}^\ell \frac{\p f^i}{\p
z^\alpha}\frac{\p f^j}{\p \bar z^\beta}\frac{\p f^k}{\p
z^\gamma}dz^\gamma\ts e_\ell\ee where the last step follows by
Bianchi identity.
 Therefore $$\left\la \sq
[R_{1,1}^{f^*(TN)},\Lambda]\p f, \p f\right\ra
=-2R_{ijk\ell}\left(h^{\alpha\bar\beta}\frac{\p f^i}{\p
z^\alpha}\frac{\p f^\ell}{\p \bar
z^\beta}\right)\left(h^{\gamma\bar\delta}\frac{\p f^j}{\p
z^\gamma}\frac{\p f^k}{\p \bar z^\delta}\right). $$ \eproof
 \elemma

\btheorem[{\cite[Sampson]{Sa2}}] Let $f:(M,h)\>(N,g)$ be a harmonic
map from a compact K\"ahler manifold $(M,h)$ to a Riemannian
manifold $(N,g)$. Then $rank_\R df\leq 2$  if $(N,g)$ has strongly
Hermitian-negative curvature. \bproof By formula (\ref{bochner}) for
the vector bundle $E=f^*(TN)$ when $(M,h)$ is K\"ahler, \beq
\Delta_{\bp_E}\bp f=\Delta_{\p_E}\bp f+\sq [R^E,\Lambda]\bp f. \eeq
If $f$ is harmonic, i.e. $\bp_E^*\bp f=0$, we obtain,
$\Delta_{\bp_E}\bp f=0$. That is \beq 0=\|\p_E\bp f\|^2+\int_M
Q_0\frac{\omega_h^m}{m!} \eeq If $(N,g)$ has strongly
Hermitian-negatve curvature, i.e. $Q_0\geq 0$ pointwisely, then
$Q_0=0$. Hence we get $rank_{\R} df\leq 2$. \eproof \etheorem

Now we go back  to work on the  Hermitian (domain) manifold $(M,h)$.
As similar  as  Lemma \ref{identity}, we obtain

 \blemma\label{realidentity} Let $f:(M,h)\>(N,g)$ be a
smooth map from a compact Hermitian manifold $(M,h)$ to a Riemannian
manifold $(N,g)$. Then\beq \int_M\p\bp\{\bp f, \bp f\}
\frac{\omega_h^{m-2}}{(m-2)!}=4\int_M \left(|\p_E\bp f|^2-|Tr_\omega
\p_E\bp
f|^2\right)\frac{\omega^m_h}{m!}+\int_M4Q_0\cdot\frac{\omega_h^m}{m!}.
\label{siu15}\eeq  \elemma

\btheorem Let $(M,h)$ be a compact Hermitian manifold with
$\p\bp\omega^{m-2}_h=0$ and $(N,g)$ a  Riemannian manifold. Let
$f:(M,h)\>(N,g)$ be a Hermitian harmonic map, then $rank_{\R}df\leq
2$ if $(N,g)$ has strongly Hermitian-negative curvature. In
particular, if $\dim_{\C}M>1$, there is no Hermitian harmonic
immersion of $M$ into Riemannian  manifolds of constant negative
curvature. \bproof If $f$ is Hermitian harmonic, i.e.,
$Tr_\omega\p_E\bp f=0$, by formula (\ref{siu15}),\beq
\int_M\p\bp\{\bp f, \bp f\} \frac{\omega_h^{m-2}}{(m-2)!}=4\int_M
|\p_E\bp
f|^2\frac{\omega^m_h}{m!}+\int_M4Q_0\cdot\frac{\omega_h^m}{m!}. \eeq
From integration by parts, we obtain
 $$4\int_M |\p_E\bp
f|^2\frac{\omega^m_h}{m!}+\int_M4Q_0\cdot\frac{\omega_h^m}{m!} =0.$$
If $(N,g)$ has strongly  Hermitian-negative curvature, then $Q_0=0$
and so  $rank_{\R}df\leq 2$.
 \eproof \etheorem

\bcorollary Let $M=\S^{2p+1}\times \S^{2q+1}$($p+q\geq 1$) be the
Calabi-Eckmann manifold. Then there is no Hermitian harmonic
immersion of $M$ into manifolds of constant negative curvature.
\bproof By a result of Matsuo(\cite{Matsuo}), every Calabi-Eckmann
manifold has a Hermitian metric $\omega$ with $\p\bp\omega^{n-2}=0$.
\eproof \ecorollary \bremark \bd\item If $M$ is  K\"ahler,  a
Hermitian harmonic immersion is also minimal.
\item By Proposition \ref{balanced}, if the manifold $(M,h)$ is
balanced, then Hermitian harmonic map is harmonic. However, if
$\omega_h$ is balanced( i.e. $d^*\omega_h=0$) and also
$\p\bp\omega_h^{m-2}=0$, then $\omega_h$ must be
K\"ahler(\cite{MT}).

\ed \eremark

\section{Rigidity of pluri-harmonic maps }\label{pl}

\subsection{Pluri-harmonic maps from complex manifolds to K\"ahler
manifolds}\label{hdddd} Let $f:M\>(N,g)$ be a pluri-harmonic map
from the compact complex manifold $M$ to the compact K\"ahler
manifold $(N,g)$. From the definition formula (\ref{pluriharmonic}),
the pluri-harmonicity of $f$ is independent of the background metric
on the domain manifold $M$ and so we do not impose a metric there.

When the domain manifold $M$ is K\"ahler, there is a number of
results on the complex analyticity and rigidity of the
pluri-harmonic $f$, mainly due to Ohinta, Udagawa and also
Burns-Burstall-Barttolomeis (e.g.,
\cite{OU},\cite{OY},\cite{BBDR},\cite{U1},\cite{U2} and the
references therein).  The common feature in their results is that
they need even more properties of the K\"ahler manifold $M$, for
example, $c_1(M)>0$, or $b_2(M)=1$.

 Now we present our main results in this section.
\btheorem\label{pluri1} Let $f:M\>(N,g)$ be a pluri-harmonic map
from a compact complex manifold $M$ to a compact K\"ahler manifold
$(N,g)$. Then it is holomorphic or anti-holomorphic if $(N,g)$ has
non-degenerate curvature and $rank_{\R} df\geq 4$. In particular,
when $(N,g)$ has strongly negative curvature (in the sense of Siu)
and $rank_{\R} df\geq 4$, then  $f$ is holomorphic or
anti-holomorphic. \etheorem

 From the proof, we can see that this
theorem also holds if the target manifold $N$ is a compact quotient
of a bounded symmetric domain and $f$ is a submersion. 

\bproof We  fix an arbitrary Hermitian metric $h$ on $M$. Let
$E=f^*(T^{1,0}N)$ and $R^E$ be the $(1,1)$-part of the curvature
tensor of $E$. If $f$ is pluri-harmonic, i.e. $\p_E\bp f=0$, then by
the Bochner formula (\ref{bochner}), the equation
$$\Delta_{\bp_E}\bp f=\Delta_{\p_E}\bp f+\sq[R^E,\Lambda](\bp f)+(\tau^*\p_E+\p_E\tau^*)(\bp f)-(\bar\tau^*\bp_E+\bp_E\bar\tau^*)(\bp f)$$
is equivalent to \beq \bp_E\bp_E^*\bp f= \sq[R^E,\Lambda](\bp
f)-(\bp_E\bar\tau^*)(\bp f).\label{key11}\eeq On the other hand, by
Lemma \ref{realbochner}, we have the relation $
[\Lambda,\p_E]=\sqrt{-1}(\bp_E^*+\bar{\tau}^{*})$,  and so \beq
\bp_E(\bp_E^*+\bar\tau^*)\bp f=-\sq \bp_E\Lambda\p_E\bp f=0 \eeq
since $f$ is pluri-harmonic. By (\ref{key11}), we get the identity $
Q=\la\sq[R^E,\Lambda](\bp f),\bp f\ra=0. $ (Note that we get $Q=0$
without using the curvature property of $(N,g)$, which is different
from the proofs in \cite[Siu]{Siu} and \cite[Jost-Yau]{JY}!) By
formula (\ref{curvature}) and the assumption that $(N,g)$ has
non-degenerate curvature, we obtain $\p f^i\wedge \p \bar f^j=0$ for
any $i$ and $j$.
 If $rank_\R(df)\geq 4$,  by
Siu's argument (\cite{Siu}), $f$ is holomorphic or
anti-holomorphic.\eproof

By using Theorem \ref{pluri1}, we can generalize a number of results
in \cite{OU}, \cite{U1} and \cite{U2} to complex (domain) manifolds.

\bproposition\label{constant}Let $M$ be an arbitrary $m$-dimensional
$(m\geq 2)$ compact complex manifold, $(N,g)$ a compact K\"ahler
manifold and $f:M\>(N,g)$ a pluri-harmonic map. Suppose  $M$ has one
of the following properties\bd
\item $\dim_{\C} H^2(M)=0$; or
\item $\dim_{\C} H^{1,1}(M)=0$; or

\item $\dim_{\C}H^2(M)=1$ and $H^2(M)$ has a generator $[\eta]$ with
$\int_M\eta^m\neq 0$; or

\item $\dim_{\C} H^{1,1}(M)=1$ and $H^{1,1}(M)$ has a generator $[\eta]$ with
$\int_M\eta^m\neq 0$. \ed
 Then \bd

 \item  $f$ is constant if $rank_\R df<2m$. In particular, if $m>n$,
 then $f$ is constant.

\item $f$ is holomorphic or
anti-holomorphic if $(N,g)$ has non-degenerate curvature. Here, we
have no rank restriction on $df$.

\ed

\bproof If $rank_{\R}df<2m$, we can consider the following real
$(1,1)$ form \beq \omega_0=\frac{\sq}{2}g_{i\bar j}\frac{\p f^i}{\p
z^\alpha}\frac{\p \bar f^j}{\p \bar z^\beta} dz^\alpha\wedge d\bar
z^\beta=\frac{\sq}{2}g_{i\bar j} \p f^i\wedge \bp \bar f^j. \eeq If
$f$ is pluri-harmonic, by Lemma \ref{key111},
$\p\omega_0=\bp\omega_0=0=d\omega_0.$ On the other hand, when
$rank_{\R}df< 2m$, $\omega_0^m=0$.  If conditions $(1)$ or $(3)$
holds, we obtain $\omega_0=d\gamma_0$. If conditions $(2)$ or $(4)$
holds, we have $\omega=\bp\gamma_1$. In any case, by Stokes'
Theorem, $\int_C\omega_0=0$ on any closed curve $C$ of $M$. But
$\omega_0$ is a nonnegative $(1,1)$ form on $M$, we obtain
$\omega_0=0$. Therefore $\p f=0$. Similarly, by using
$$\omega_1=\frac{\sq}{2}g_{i\bar j}\frac{\p \bar f^j}{\p z^\alpha}\frac{\p f^i}{\p \bar z^\beta}dz^\alpha\wedge
d\bar z^\beta= \frac{\sq}{2} g_{i\bar j}\p\bar f^i\wedge \bp f^j$$
we know $\bp f=0$. Hence $f$ is constant. In particular, if $m>n$,
i.e.  $rank_{\R}df<2m$,  $f$ is constant.

Suppose $(N,g)$ has non-degenerate curvature. If $rank_{\R} df \geq
2m\geq 4 $, by Theorem \ref{pluri1},  then $f$ is holomorphic or
anti-holomorphic. If $rank_{\R} df < 2m$, by the proof above, we see
$f$ is constant.
 \eproof \eproposition

%
%

 \bcorollary Any pluri-harmonic map from the Calabi-Eckmann manifold
$\S^{2p+1}\times \S^{2q+1}$ to the $n$-dimensional complex space
form $N(c)$ is constant if $p+q\geq n$.  \ecorollary

The following result is well-known, (e.g.\cite{BBDR}, \cite{OY},
\cite{OU}).

\bcorollary Every pluri-harmonic map from $\P^m$ to $\P^n$ is
constant if $m>n$. \ecorollary

\vskip 1\baselineskip

\subsection{Pluri-harmonic maps from Hermitian manifolds to Riemannian
manifolds} In this subsection, we shall use similar ideas as in
Section \ref{hdddd} to study the rigidity of pluri-harmonic maps
from  Hermitian manifolds to Riemannian manifolds.

\btheorem\label{pluri2} Let $f:M\>(N,g)$ be a pluri-harmonic map
from a compact complex manifold $M$ to a  Riemannian manifold
$(N,g)$. If $(N,g)$ has non-degenerate Hermitian curvature at some
point $p$, then $rank_{\R} df(p)\leq 2$.\bproof We fix an arbitrary
Hermitian metric $h$ on $M$. If $f$ is pluri-harmonic, i.e. $\p_E\bp
f=0$, then by the Bochner formula (\ref{bochner}), the equation
$$\Delta_{\bp_E}\bp f=\Delta_{\p_E}\bp f+\sq[R^E,\Lambda](\bp f)+(\tau^*\p_E+\p_E\tau^*)(\bp f)-(\bar\tau^*\bp_E+\bp_E\bar\tau^*)(\bp f)$$
is equivalent to $ \bp_E\bp_E^*\bp f= \sq[R^E,\Lambda](\bp
f)-(\bp_E\bar\tau^*)(\bp f).$ By a similar argument as in Theorem
\ref{pluri1}, we obtain
 $ Q_0=\la\sq[R^E,\Lambda](\bp f),\bp
f\ra=0.$  That is \beq R_{ijk\ell}\left(h^{\alpha\bar\beta}\frac{\p
f^i}{\p z^\alpha}\frac{\p f^\ell}{\p \bar
z^\beta}\right)\left(h^{\gamma\bar\delta}\frac{\p f^j}{\p
z^\gamma}\frac{\p f^k}{\p \bar z^\delta}\right)=0. \eeq If the
curvature tensor $R$ of $(N,g)$ is non-degenerate at some point
$p\in M$, then  the complex rank of the matrix $
\left(h^{\alpha\bar\beta}\frac{\p f^i}{\p z^\alpha}\frac{\p f^j}{\p
\bar z^\beta}\right)$
 is $\leq 1$, i.e. $rank_{\R}df(p)\leq 2$.
 \eproof
\etheorem

\bproposition \label{constant2} Let $M$ be an arbitrary
$m$-dimensional($m\geq 2$) compact complex manifold, $(N,g)$ a
Riemannian manifold and $f:M\>(N,g)$  a pluri-harmonic map. Suppose
$M$ has one of the following properties \bd
\item $\dim_{\C} H^2(M)=0$; or \item $\dim_{\C} H^{1,1}(M)=0$; or

\item $\dim_{\C}H^2(M)=1$ and $H^2(M)$ has a generator $[\eta]$ with
$\int_M\eta^m\neq 0$; or

\item $\dim_{\C} H^{1,1}(M)=1$ and $H^{1,1}(M)$ has a generator $[\eta]$ with
$\int_M\eta^m\neq 0$; \ed

 then \bd \item $f$ is constant if $rank_\R df<2m$. In particular, if $m>n$,
 then $f$ is constant.

\item $f$ is  constant if $(N,g)$ has non-degenerate curvature.

\ed

  \bproof
Assume $rank_{\R} df<2m$. We can consider
 $
\omega_0=\frac{\sq}{2}g_{i\bar j} \p f^i\wedge \bp  f^j.$ By a
similar proof as in Proposition \ref{constant}, we obtain $\bp f=0$,
and so $f$ is a constant. On the other hand, if $(N,g)$
 has non-degenerate curvature, then $rank_\R df\leq 2<2m$, hence $f$
 is constant.
\eproof \eproposition

 \bcorollary \bd\item Any pluri-harmonic map from the Calabi-Eckmann manifold
$\S^{2p+1}\times \S^{2q+1}$ to the real space form $N(c)$ is
constant if $p+q\geq 1$.

\item Any pluri-harmonic map from  $\C\P^n$ to the real space form $N(c)$ is
constant if $n\geq 2$.

 \ed\ecorollary

\end{document}